\documentclass[12pt,oneside,reqno]{amsart}
\usepackage{graphicx}
\usepackage{caption,subcaption}
\usepackage{dcolumn}
\usepackage{amsmath,amsthm,amssymb}
\usepackage{todonotes}
\usepackage{url}
\usepackage{placeins}

\usepackage{setspace}
\usepackage{array,booktabs}

\renewcommand{\Re}{\operatorname{Re}}
\renewcommand{\Im}{\operatorname{Im}}

\def\be{\begin{equation}}
\def\ee{\end{equation}}

\newtheorem*{theorem*}{Theorem}
\newtheorem*{conj*}{Conjecture}
\theoremstyle{remark}
\newtheorem*{rem*}{Remark}

\begin{document}

\title[Sine-Gordon on a wormhole]{Sine-Gordon on a wormhole}

\author{Piotr Bizo\'n}
\address{Institute of Theoretical Physics, Jagiellonian
University, Krak\'ow}
\email{bizon@th.if.uj.edu.pl}
\author{Maciej Dunajski}
\address{DAMTP, University of Cambridge, Cambridge}
\email{m.dunajski@damtp.cam.ac.uk}

\author{Micha\l{} Kahl}
\address{Institute of Theoretical Physics, Jagiellonian
University, Krak\'ow}
\email{michal.kahl@alumni.uj.edu.pl}

\author{Micha\l{} Kowalczyk}
\address{Departamento de Ingenier\'ia Matem\'atica and Centro de Modelamiento Matem\'atico\\
Universidad de Chile\\  Santiago, Chile}
\email{kowalczy@dim.uchile.cl}

\date{\today}%
\begin{abstract}
In an attempt to understand the soliton resolution conjecture, we consider the Sine-Gordon equation on a spherically symmetric wormhole spacetime. We show that within each topological sector (indexed by a positive integer degree $n$) there exists a unique linearly stable soliton, which we call  the $n$-kink.  We give numerical evidence that the $n$-kink is a global attractor in the evolution of any smooth, finite energy solutions of degree $n$. When the radius of the wormhole throat $a$ is large enough, the  convergence to the $n$-kink is shown to be governed  by internal modes that slowly decay  due to the resonant transfer of energy to radiation.  We compute the exact asymptotics of this relaxation process for the $1$-kink using the Soffer-Weinstein weakly nonlinear perturbation theory.
\end{abstract}
\maketitle

\section{Introduction}

If a solution of an evolution equation exists for all times~$t$, then it is natural to ask how it behaves as $t\rightarrow \infty$. This question is particularly interesting if the equation admits solitons (spatially localized, finite energy solutions) because they may appear as late-time attractors. For nonlinear dispersive wave equations it is believed that for any reasonable (e.g. smooth and finite energy) generic initial data, the solution eventually resolves into a superposition of a radiative component plus a finite number of solitons. This belief, known as the soliton resolution conjecture \cite{avy}, is fairly well understood for small perturbations of solitons \cite{tao, KMM}, however little is known  in the non-perturbative regime (but the one-dimensional completely integrable equations where  solutions can be computed explicitly via inverse scattering methods \cite{ES, Sch} and few results in higher dimensions, e.g. \cite{dkm, klls}).

 In an attempt to understand the soliton resolution conjecture in a simple setting,  two of us proposed in \cite{BK} to study nonlinear waves propagating on a spherically symmetric curved spacetime with the metric
\begin{equation}\label{worm}
  ds^2=-dt^2+ dr^2 + (r^2+a^2) d\omega^2 \,,
\end{equation}
where $(t,r)\in \mathbb{R}^2$, $d\omega^2$ is the round metric on the unit two-sphere, and $a$ is a positive constant. This spacetime, introduced by Ellis \cite{ellis} and Bronnikov \cite{bron}, is the simplest example of a wormhole geometry that has two asymptotically flat ends at $r \rightarrow \pm \infty$ connected by a spherical throat (minimal surface)  of area $4\pi a^2$ at $r=0$.

\begin{figure}
\includegraphics[width=.5\textwidth]{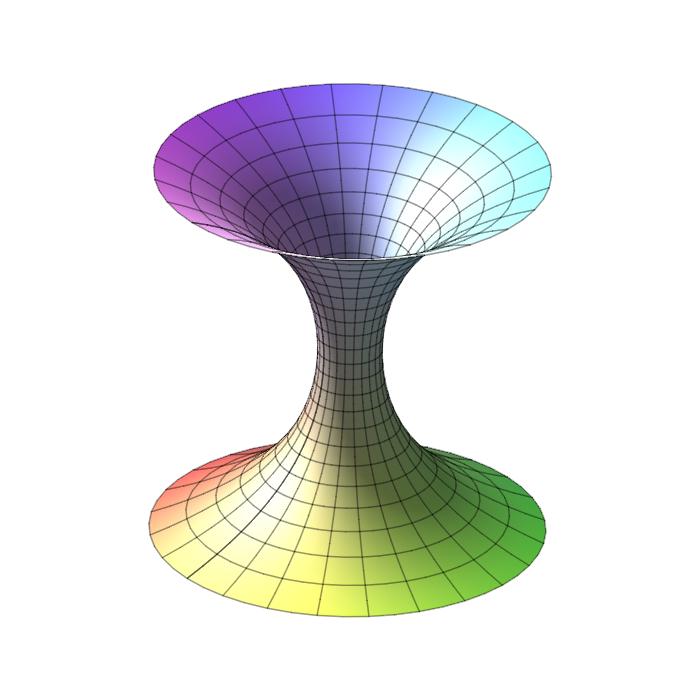}
\includegraphics[width=.48\textwidth]{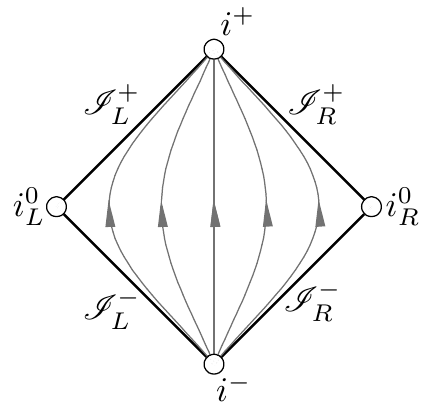}
\caption{Left panel: the isometric embedding of the constant time  equatorial cross--section of the wormhole in $\mathbb{R}^3$. Right panel: the conformal diagram of the wormhole spacetime. The boundaries of each side of  the diamond, denoted by $\mathcal{J}^{\pm}_{L,R}$, represent future/past ($t\rightarrow \pm \infty$) and left/right ($r\rightarrow \mp \infty$) null infinities. }
\end{figure}
Although the wormhole belongs more to science-fiction  than physics \cite{MT}, it has  a number of features that makes it an attractive testing ground for  the soliton resolution conjecture.
 First, there is no singularity at $r=0$ which basically ensures global well-posedness for dispersive equations with coercive nonlinearities. Second, due to  the presence of the length scale~$a$, Derrick's non-existence scaling arguments are evaded and solitons do exist (often in abundance) for some nonlinearities. Third, if a soliton exists, it is completely rigid so no modulation analysis is needed. Finally,
  the equations posed on the wormhole combine the simplicity of  one-dimensional equations on the whole real line with the three-dimensional dispersive properties.

  In \cite{BK} the soliton resolution conjecture was formulated and verified numerically for equivariant wave maps from the wormhole \eqref{worm} into the 3-sphere. In addition,  the  rate of convergence to the soliton (which in this case is a harmonic map from a $t=\text{const}$ hypersurface of the wormhole into the 3-sphere) was computed by perturbation methods. Subsequently, the conjecture made in \cite{BK} was proved (without a decay rate, though) by Rodriguez \cite{R1,R2} via the concentration-compactness method as in \cite{klls}. The key ingredient in getting these results was the fact that the linearized perturbations around the solitons decay in time.   In this paper, we consider a different nonlinearity for which the latter property does not hold and the asymptotic stability of solitons is an inherently nonlinear phenomenon.

On the wormhole spacetime we consider a real scalar field $\phi$ obeying the  semilinear wave equation\footnote{For dimensional reasons the nonlinear term must have the form $\ell^{-2} \sin(2\phi)$, where $\ell$ is a fixed scale of length. Hereafter, we set $\ell=1$ by the choice of the unit of length.}
\begin{equation}\label{sgeq}
  \Box_g \phi + \sin(2\phi)=0,
\end{equation}
where $\Box_g$ is the wave operator associated with the metric \eqref{worm}.
Assuming that $\phi=\phi(t,r)$, we get
\begin{equation}\label{eq}
  \ddot \phi=\phi''+\frac{2r}{r^2+a^2}\, \phi' - \sin(2\phi),
\end{equation}
where an overdot and prime denote derivatives with respect to $t$ and $r$. For $a=\infty$ this equation reduces to the one-dimensional Sine-Gordon equation
\begin{equation}\label{eqsn}
  \ddot \phi=\phi'' - \sin(2\phi),
\end{equation}
 which is completely integrable, hence for large values of $a$ equation \eqref{eq} can be viewed as a non-integrable perturbation of \eqref{eqsn}\footnote{The $\phi^4$ model on the wormhole was considered in \cite{aw} following an earlier unpublished version of our paper. The key difference which makes the $\phi^4$ model less interesting than Sine-Gordon  is that the former is non-integrable already in flat space. Equation \eqref{eq} is indeed not integrable for any finite $a$ as it does not possess the Painlev\'e property.}.

The conserved energy associated with equation \eqref{eq} reads
\begin{equation}\label{energy}
  E =  \int\limits_{-\infty}^{\infty} \left(\frac12 \dot\phi^2 + \frac12 \phi'^2 +\sin^2{\!\phi} \right)(r^2+a^2)  \, dr\,.
\end{equation}
Finiteness of energy requires that $\phi(t,-\infty)=n_-\pi$, $\phi(t,\infty)=n_+\pi$, where $n_-$ and $n_+$ are integers. Without loss of generality we choose $n_-=0$; then $n=n_+$ determines the topological degree of the solution (which is preserved in the evolution).

The goal of this paper is to describe the asymptotic behaviour of solutions of equation \eqref{eq} for $t\rightarrow \infty$.
 Due to the dissipation of energy by dispersion,  solutions are expected to settle down to stationary states, in accord with the soliton resolution conjecture.
 In section~2 we prove that for each degree $n$ there exists a unique smooth, finite-energy stationary solution, which we call the $n$-kink.  The linear stability of the $n$-kinks is analyzed in section~3. We show that the spectrum of the linearized operator around the $n$-kink has no negative or zero eigenvalues, hence the $n$-kink is linearly stable.  However, for sufficiently large $a$ there are $n$ positive eigenvalues in the mass gap between zero and the bottom of the continuous spectrum. These positive eigenvalues  give rise to internal  modes that oscillate harmonically and therefore prevent asymptotic stability of kinks at the linear level (if $a$ is large enough). Nonetheless, it is expected that the $n$-kink is asymptotically stable thanks to the nonlinear resonant damping of internal modes, as described by Soffer and Weinstein in \cite{SW}. In section~4 we use their weakly nonlinear perturbation method to derive the decay rate of the internal mode for the $1$-kink. Finally, in section~5 we give numerical evidence for the soliton resolution conjecture and verify the predictions of perturbative computations.  As in \cite{BK}, we solve equation \eqref{eq}  numerically using  the hyperboloidal formulation of the initial value problem.   This approach allows us to reach very long times of evolution in a reasonable computational time. In the appendix we give some details of the computation of  parameters of kinks.

 \section{Kinks}
 Time-independent solutions $\phi=\phi(r)$ of equation  \eqref{eq} satisfy the ordinary differential equation
\begin{equation}\label{ode}
  \phi''+\frac{2r}{r^2+a^2}\, \phi'- \sin(2\phi) =0\,.
\end{equation}
This equation can be viewed  as the equation of motion, with `time' $r$, for the unit mass particle moving in the potential $-\sin^2{\!\phi}$ and subject to friction with the time-dependent friction coefficient $\frac{2r}{r^2+a^2}$.
The solution of degree $n$ corresponds to the trajectory whose projection on the phase plane $(\phi,\phi')$ starts from the saddle point $(0,0)$ at $r=-\infty$ and goes to the saddle point $(n\pi,0)$ for $r=+\infty$. The existence and uniqueness  of such a connecting trajectory for each $n$ follows from an elementary shooting argument. For example, let the particle be located at $\phi=\pi/2$ for $r=0$. If the velocity $b=\phi'(0)$ is too small, then the particle will never reach the hilltop at $\phi=\pi$, while if $b$ is sufficiently large it will roll over the hilltop. By continuity, there must be a critical velocity $b_1$ for which the particle reaches the hilltop in infinite time (obviously, by the uniqueness of trajectories, the particle cannot stop at the hilltop in finite time). Due to reflection symmetry $r\rightarrow -r$, the particle sent backwards in time reaches $\phi=0$ for $r\rightarrow-\infty$, giving the desired connecting trajectory with $n=1$\footnote{If $0<\phi(0)<\pi/2$, then by the same shooting argument there exists a velocity $\phi'(0)$ such that the particle tends to $\pi$ as $r\rightarrow \infty$, however going backwards this particle will overshoot $0$ and end up at $-\pi/2$ for $r\rightarrow-\infty$. Thus, there are no asymmetric kinks.}.  Repeating this argument for higher $n$ we get a countable family of unique connecting trajectories $\phi_n(r)$ which are symmetric with respect to the midpoint $\phi(0)=n\pi/2$, that is $\phi_n(r)+\phi_n(-r)=n\pi$. Near $r=0$
\begin{equation}\label{expans}
  \phi_n(r)=\frac{n\pi}{2} + b_n r +\mathcal{O}(r^3)\,,
\end{equation}
where the parameter $b_n$ uniquely determines the trajectory. For $r\rightarrow\mp\infty$ the leading asymptotics are, respectively
\begin{equation}\label{expans_inf}
  \phi_n(r)\sim  -  \frac{c_n}{r}\,e^{\sqrt{2} r} \quad\mbox{and}\quad \phi_n(r)\sim n\pi -  \frac{c_n}{r}\,e^{-\sqrt{2} r}  \,,
\end{equation}
where the parameter $c_n$ is determined by $b_n$.
In the following, we shall refer to the stationary solution $\phi_n$ as the $n$-kink. Fig.~2 depicts sample profiles of $n$-kinks for $n=1,2$.
 A few values of the parameters $b_n$ and $c_n$ are listed in
Table~I  for different values of $a$ and $n$. While the numerical computation of the parameters $b_n$ is straightforward by means of the shooting method, the computation of the parameters $c_n$ is more difficult because the leading asymptotic behavior \eqref{expans_inf} is only the first term of the asymptotic series which has to be summed to give an accurate approximation of the solution. The details of this computation are given in the appendix.

\begin{figure}[h]
\includegraphics[width=.49\textwidth]{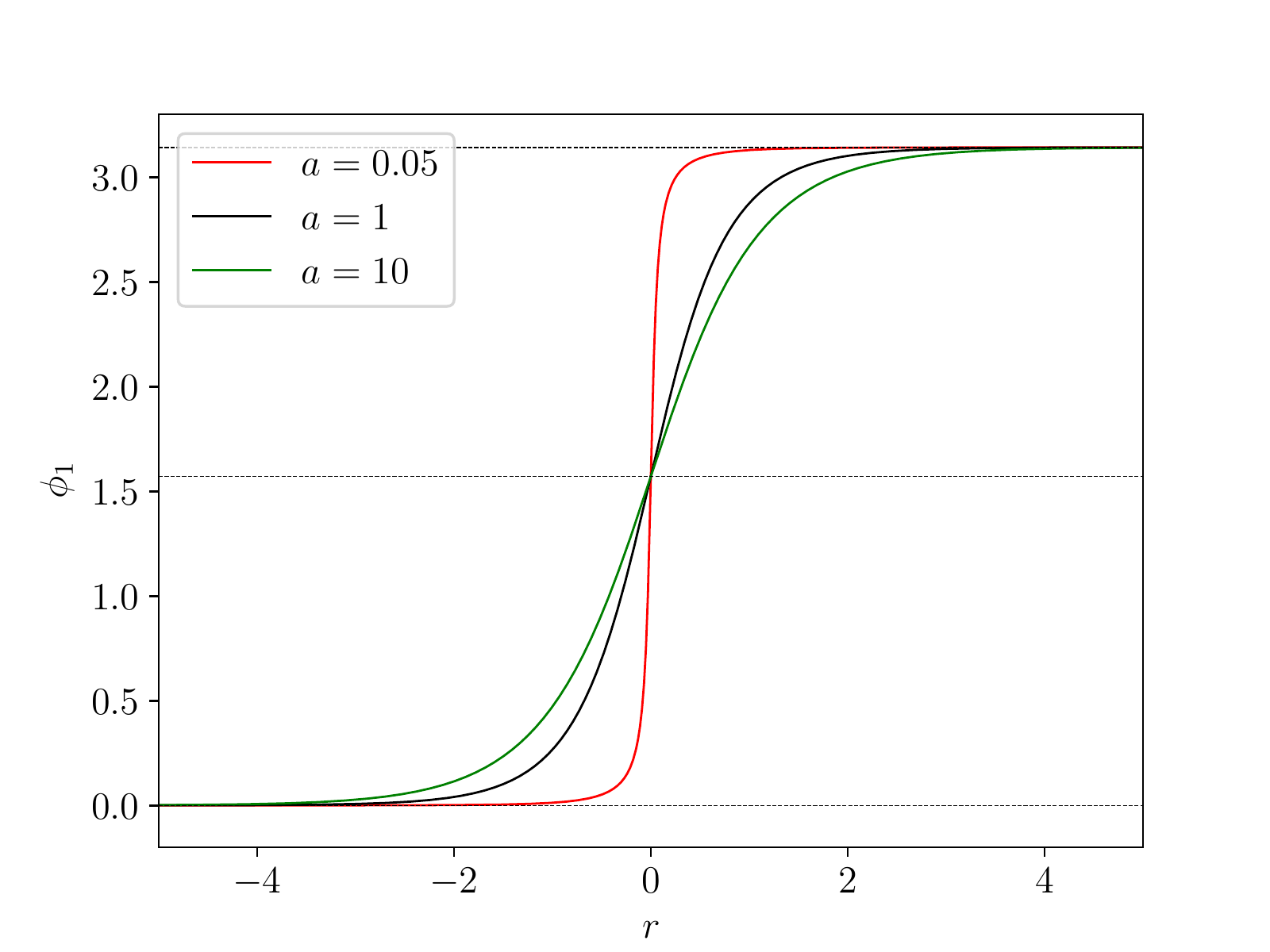}
\includegraphics[width=.49\textwidth]{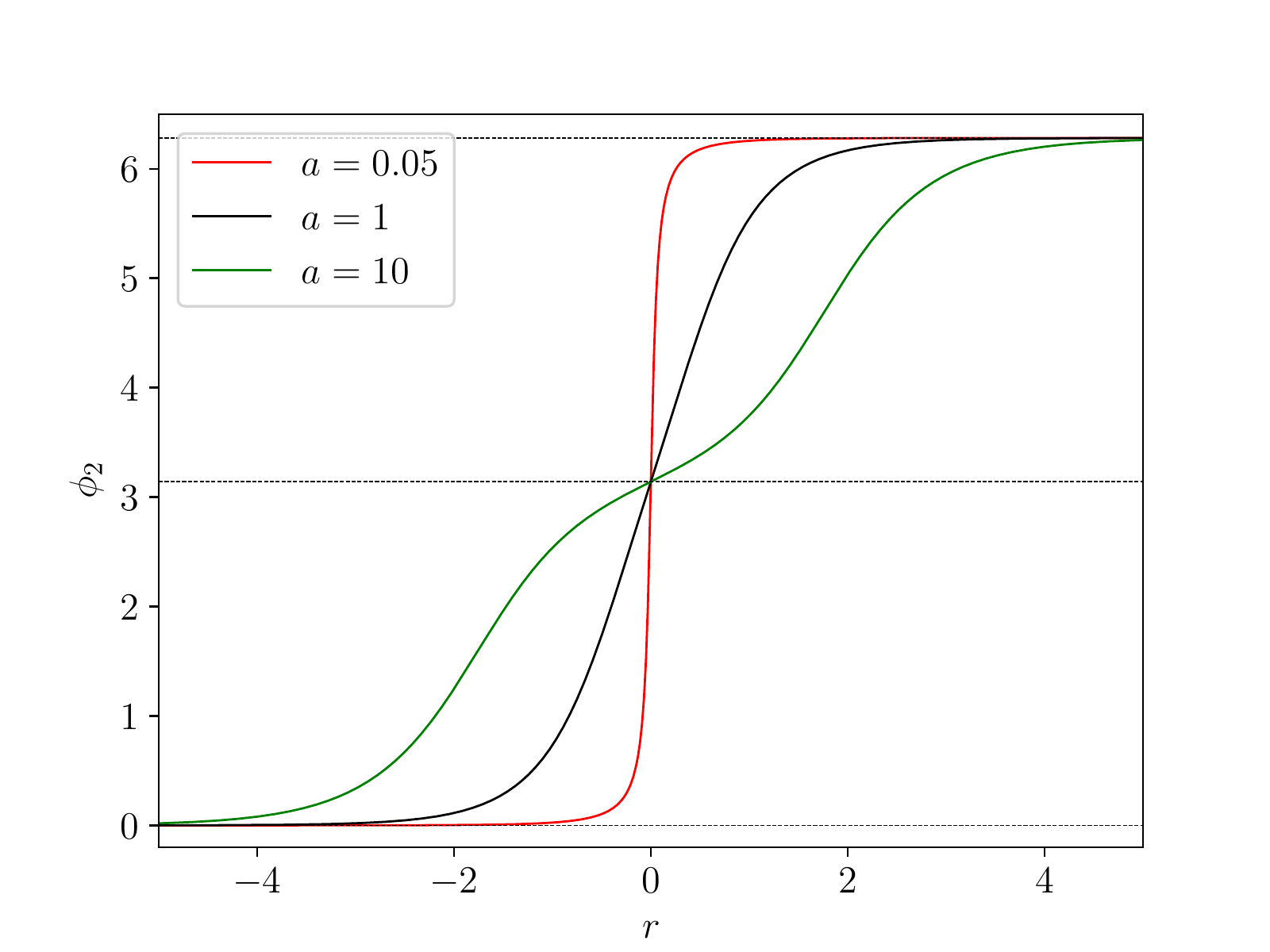}
\caption{{\small{Profiles of the $n$-kinks for $n=1,2$.}}}
\end{figure}
\begin{table} [h]
   \begin{tabular}{|c|ccc|}
    \toprule
    $a $ & $(b_{1},c_{1})$ & $(b_{2},c_{2})$ & $(b_{3},c_{3})$ \\
    \midrule
    1 & (2.0163,\,1.5054) & (2.8709,\,4.2523) & (4.3285,\,8.5162) \\
    2 & (1.6152,\,3.4063) & (1.6531,\,13.109) & (2.7121,\,33.218) \\
    3 & (1.5123,\,5.3885) & (1.1993,\,26.592) & (2.1862,\,82.056) \\
    \bottomrule
  \end{tabular}
  \vskip 2ex
  \caption{{\small{Parameters $b_{n}$ and $c_n$ for three values of $a$.}}}
  \label{tab:1}
\end{table}

Next, we shall derive analytic approximations of the kink solutions for large and small values of $a$.
These approximations will be used below in the stability analysis of the kink.
\vskip 0.2cm
 {\bf Large $a$ approximation:} For $a=\infty$, equation \eqref{ode} reduces to the static Sine-Gordon equation
\begin{equation}\label{ode-sg}
  \phi''- \sin(2\phi) =0\,,
\end{equation}
whose unique (modulo translation) soliton solution is the Sine-Gordon kink
\be
\label{sg-kink}
H(r)=2 \arctan{e^{\sqrt{2} r}}\,.
\ee
Let $\varepsilon=1/a\ll1$ and write
\be
\label{eq-psi1}
\phi(r)=H(r)+\varepsilon^2 \psi(r)+\mathcal{O}\left(\varepsilon^4\right).
\ee
Substituting this expansion into equation \eqref{ode} and collecting terms of order $\varepsilon^2$, we obtain
\be
\label{eq-psi}
{\psi}''-2\cos{(2H)}\psi=-2r H'\,.
\ee
 The solution $H'$ of the homogeneous equation is even, while the right hand side is odd, hence the Fredholm solvability condition is satisfied and consequently there is a unique  solution that is odd and  decays at $\pm \infty$. This solution can be written in closed form (using a polylogarithmic function) but we refrain from displaying it here because it will not be used  below.
We confirmed numerically that the approximation  $\phi_1(r)\approx H(r)+\psi(r)/a^2$ is very accurate for sufficiently large $a^2$ (and $r^2<a^2$).

The $n$-kink for large $a$ can be approximated by a superposition of $n$ well-separated Sine-Gordon kinks. For example, for $n=2$ we have
\be
\label{2_kink}
\phi_2(r) \approx H(r-R)+H(r+R).
\ee
The dependence of the separation parameter $R$ on $a$ can be calculated as follows. Multiplying equation \eqref{ode} by $\phi'$ and integrating from $r=0$ to $r=\infty$, we get
\be
\label{integral_a}
\frac{1}{2}\phi'(0)^2-\sin^2{(\phi(0))}={\int_0}^{\infty}
 \frac{2r}{r^2+a^2}{\phi'}^2dr.
\ee
In terms of the mechanical analogy this equation represents the balance between the initial  energy of the fictitious particle and the energy lost by friction.
 Substituting (\ref{2_kink}) into (\ref{integral_a}) and
assuming that $R$ and $a$ are large, we get at the leading order
$
 16 e^{-2\sqrt{2} R}=\gamma/a^2,
$
where $\gamma$ is a constant. Thus,
\be
\label{logR}
R \approx \frac{1}{\sqrt{2}}\ln{(a)} \quad \mbox{for} \quad a\gg 1.
\ee
Similar large-$a$ approximations can be given for $n$-kinks with larger $n$.
\vskip 0.2cm
{\bf Small $a$ approximation:}
Changing variables to $\rho=r/a$ and $g(\rho)=\phi(r)$, and taking the limit $a\rightarrow 0$, we get the linear equation
\be
\label{scalar_lap}
\frac{d^2 g}{d\rho^2}+\frac{2\rho}{\rho^2+1}\frac{d g}{d\rho}=0,
\ee
whose two linearly independent solutions are $g_1=1$ and  $g_2=\arctan{\rho}$. Thus, for $a\ll 1$ the $n$-kink is approximated by
\be
\label{small_a}
\phi_n(r) \approx n\Big(\frac{\pi}{2}+ \arctan\left(r/a\right)\Big).
\ee

\section{Linear perturbations}
In this section we analyze  linear stability of kinks $\phi_n(r)$.
Let
\be \label{ls}
\phi(t,r)=\phi_n(r)+ (r^2+a^2)^{-1/2} u(t,r),
\ee
where the perturbation $u$ is assumed to be small.
Plugging this into equation \eqref{eq}, dropping nonlinear terms in $u$, and assuming harmonic time dependence $u(t,r)=e^{-i\omega t} v(r)$, we get the eigenvalue problem for the one-dimensional Schr\"odinger operator
\begin{equation}
\label{eq1}
L_n v\equiv \Big(-\frac{d^2}{dr^2}+2+V_n(r)\Big)v=\omega^2 v
\end{equation}
with the potential
\be
\label{potentialV}
V_n(r)=-4 \sin^2{\!\phi_n(r)}+\frac{a^2}{(r^2+a^2)^2}.
\ee

\begin{figure}[h]
\includegraphics[width=0.49\textwidth]{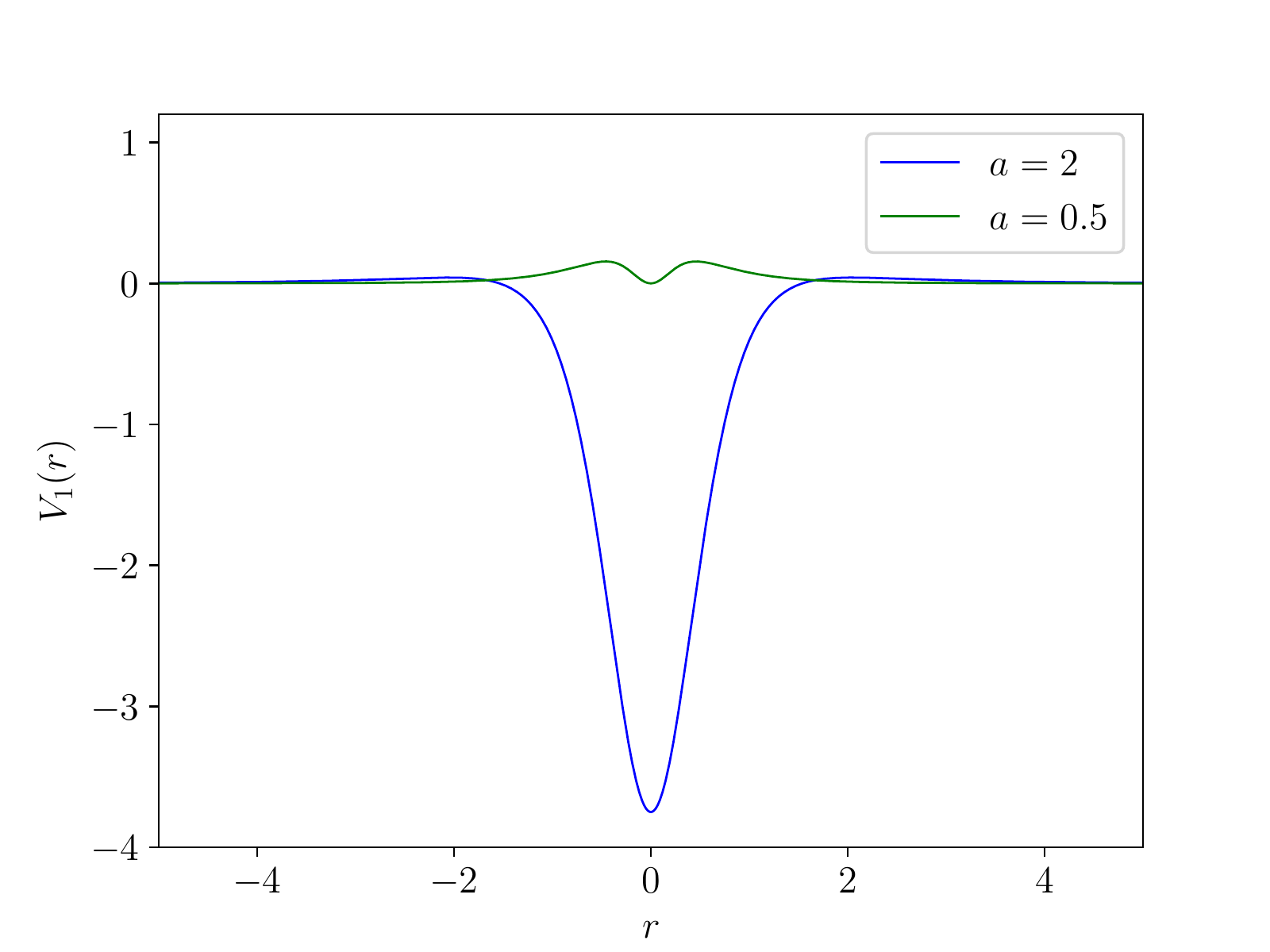}
    \includegraphics[width=0.49\textwidth]{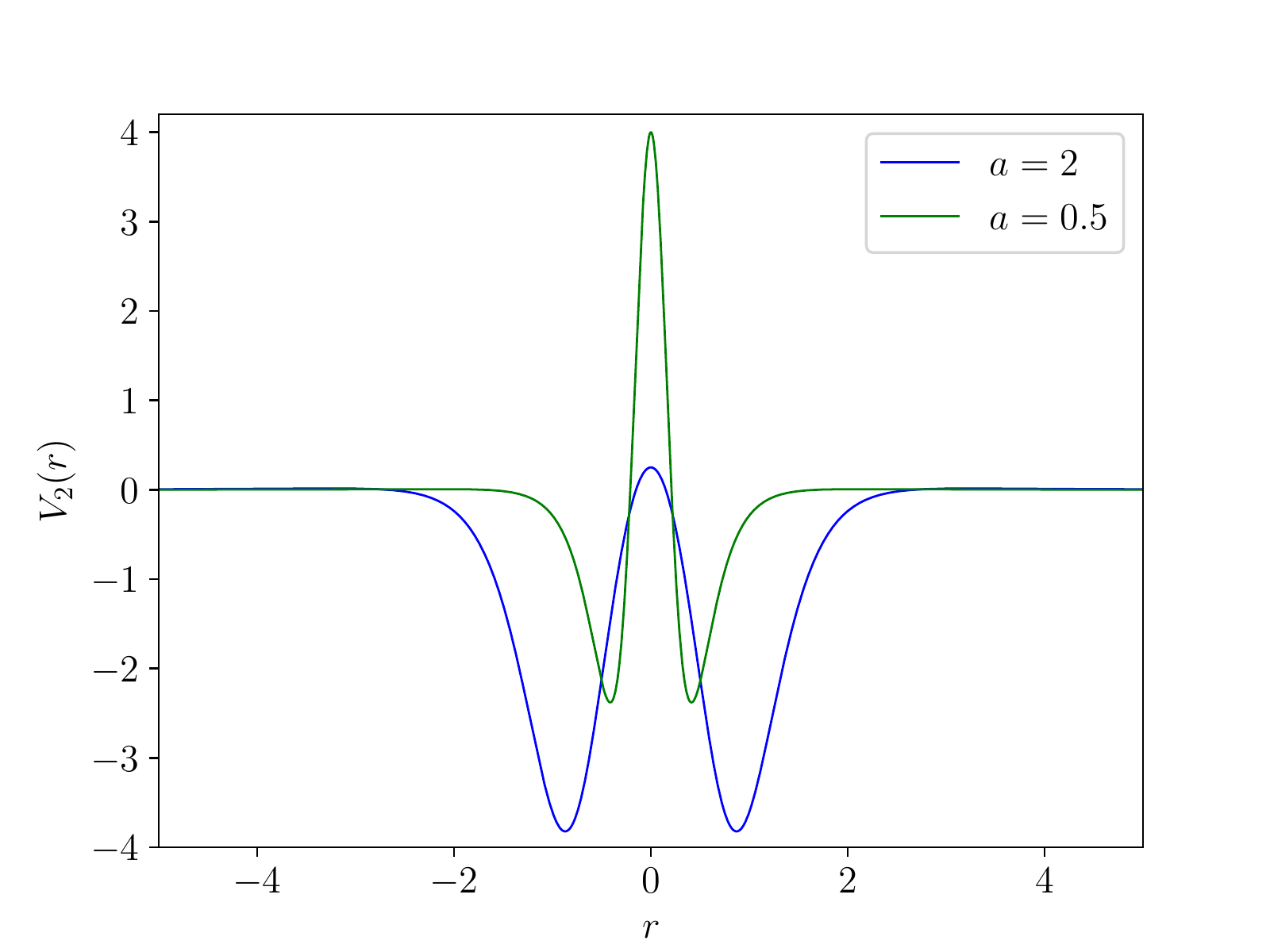}
    \caption{{\small{The potential $V_n(r)$ for $n=1,2$ and two values of $a$.}}}
    \end{figure}

Let us first consider the case $n=1$. As discussed above, in this case the 1-kink $\phi_1(r)$ tends for $a\rightarrow \infty$ to the Sine-Gordon kink $H(r)$  given by \eqref{sg-kink}
for which the corresponding linear stability operator is
\be
\label{L-sg}
L=-\frac{d^2}{dr^2}+2-4 \sin^2{\!H(r)}=-\frac{d^2}{dr^2} + 2-\frac{4}{\cosh^2(\sqrt{2} r)}\,.
\ee
This operator has a continuous spectrum $\omega^2\geq 2$ (the bottom $\omega^2=2$ is a resonance) and a single  eigenvalue $\omega^2=0$ which is due to translation symmetry of the Sine-Gordon equation;  the associated normalized  eigenfunction (zero mode) is given by
\be\label{zeromode}
v_0(r)=2^{-3/4} H'(r)=\frac{2^{-1/4}}{\cosh(\sqrt{2}r)}.
\ee
The absence of eigenvalues in the gap $(0,2)$ is believed to be intimately tied with the complete integrability of the Sine-Gordon equation \cite{kpcp}.

The operator $L_1$ can be viewed as a compact perturbation of $L$ so it has the same continuous spectrum $\omega^2\geq 2$ but the discrete spectra are different. We claim that $L_1>L$. To show this,
let us observe that $\phi_1(r)>H(r)$ for $r>0$ (and by the reflection symmetry $0<\phi_1(r)<H(r)$ for $r<0$).
This fact is evident within our mechanical analogy because the fictitious particle corresponding to $\phi_1$ is subject to friction while the one corresponding to $H$ moves without friction. To see this, note that $\phi_1(0)=H(0)=\pi/2$ and $\phi'_1(0)>H'(0)$, hence $\phi_1(r)>H(r)$ for small $r>0$. In fact, this inequality holds for all $r>0$ because the $H$-particle cannot overtake the $\phi_1$-particle (since at the overtake point the $\phi_1$-particle would have smaller kinetic energy than the $H$-particle and  could not reach the hilltop). Since $-\sin^2{\phi}$ is decreasing for $\phi\in [0,\pi/2]$ and increasing for $\phi\in [\pi/2,\pi]$, it follows that $-\sin^2{\phi_1(r)}>-\sin^2{H(r)}$ for all $r$ which by \eqref{eq1} and \eqref{L-sg} implies that $L_1>L$. Since $L$ has exactly one eigenvalue at $\omega^2=0$, it follows (see, e.g., Corollary 4.11 on page 119 in \cite{Te}) that $L_1$ has no negative eigenvalues and at most one positive eigenvalue in the gap $(0,2)$.

To obtain a more quantitative information about the gap eigenvalue of $L_1$ for large values of $a$ we seek a perturbative solution of the eigenvalue problem \eqref{eq1} in the form
\be
\label{vpert}
v=v_0+\varepsilon^2 v_1+\mathcal{O}\left(\varepsilon^4 \right),\quad \omega^2=c\, \varepsilon^2 +\mathcal{O}\left(\varepsilon^4\right),
\ee
where $\varepsilon=1/a\ll 1$. Inserting \eqref{eq-psi1} and \eqref{vpert} into \eqref{eq1}, at the zero order we obtain
$L v_0=0$, while at  the order $\varepsilon^2$ we get
\[
L v_1=c v_0+2r {v_0}'+4\sin(2H)\psi v_0.
\]
The right hand side must be orthogonal to $v_0$ which yields
\be
\label{c}
c=1-4\int_{-\infty}^{\infty}\sin{(2H(r))}\psi(r){v_0(r)}^2 dr\,.
\ee
To calculate the above  integral we differentiate equation \eqref{eq-psi}
\begin{equation}\label{eq-psid}
 L \psi' = 4 \sin(2H) H' \psi + 2 (r H')'\,.
\end{equation}
Taking the inner product with $v_0$ and integrating by parts we get
\begin{equation}\label{identity}
4 \int_{-\infty}^{\infty} \sin(2H(r)) \psi(r) v_0(r)^2 \,dr = -\int_{-\infty}^{\infty} v_0(r)^2\,dr = -1.
\end{equation}
Substituting this into \eqref{c} we obtain $c=2$.

 As $a$ decreases, the potential well gets shallower and eventually becomes a barrier, as follows from the small-$a$ approximation of kinks \eqref{small_a}; see Fig.~3. Accordingly, as shown in Fig.~4, the eigenvalue $\omega^2$ migrates through the gap $(0,2)$ and disappears into the continuous spectrum  for $a$ smaller than some critical value $a^*$ (for $a=a^*$ there is a resonance at the bottom of the continuous spectrum)\footnote{A similar behavior of the gap eigenvalue was found for some geometric wave equations on the hyperbolic space \cite{LOS}.}. Numerically, we find that $a^*\approx 0.536$.

 \begin{figure}[h]
\includegraphics[width=0.6\textwidth]{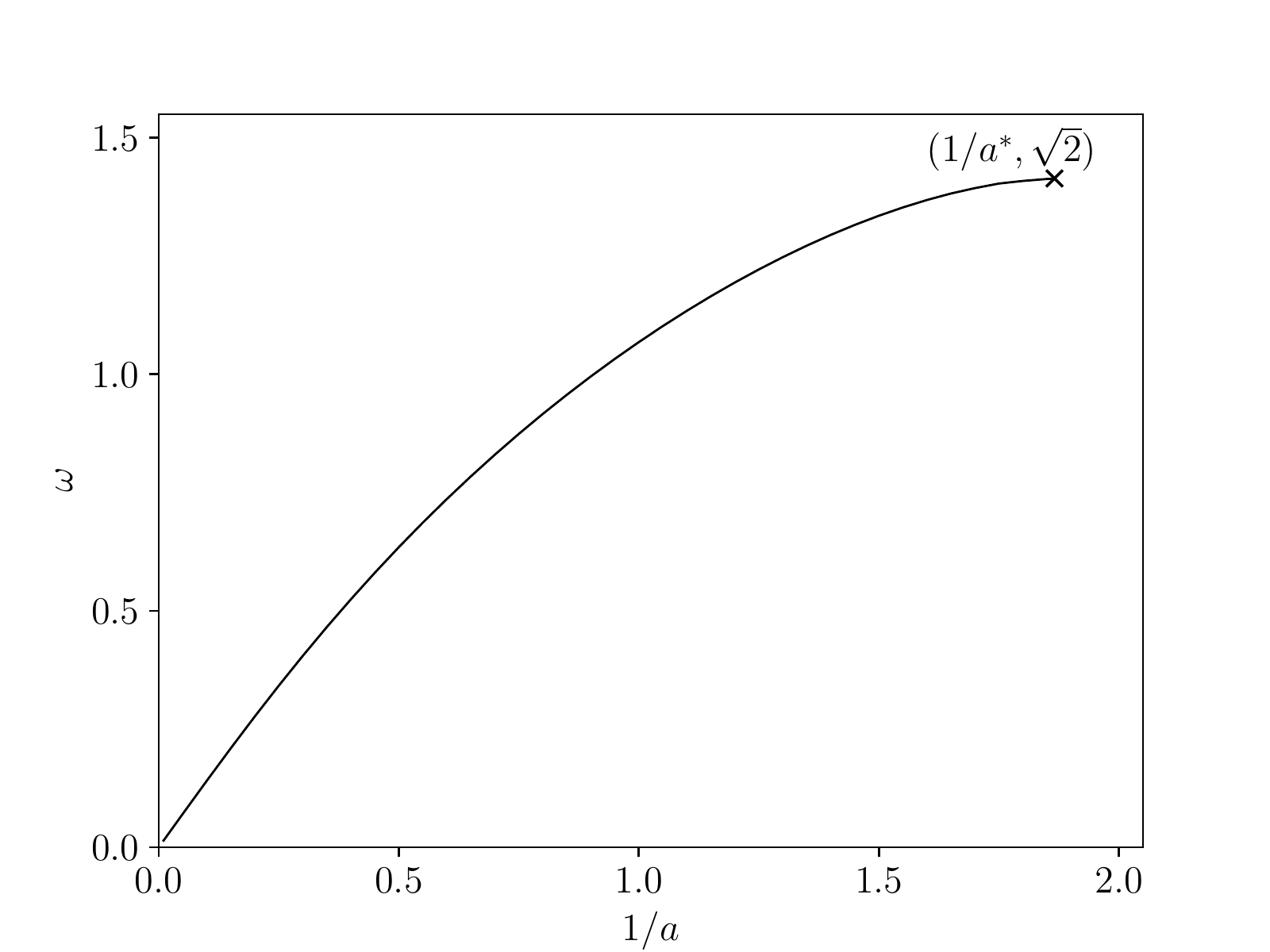}
    \caption{{\small{The frequency of the internal mode of the $1$-kink as a function of $1/a$.}}}
    \end{figure}

    For a given $n\geq 2$ and sufficiently large $a$, it follows from the large-$a$ approximation of $n$-kinks  that the potential has the form of $n$ wells equally separated by the distance $\propto \ln(a)$. Consequently, if $a$ is large enough there are  $n$ gap eigenvalues
     \[
     0< \omega_1^2 <\dots < \omega_n^2  <2.
     \]
     As $a$ decreases, the potential wells go up (see Fig.~3) and the gap eigenvalues  disappear  one by one into the continuous spectrum at certain critical values $a^*_1<a^*_2<\dots<a^*_n$.
         For example, for $n=2$  we find numerically $a^*_1\approx0.39$ and $a^*_2\approx0.81$.

\section{Weakly nonlinear dynamics near the $1$-kink}
In this section we study solutions of equation \eqref{eq} for initial data near the kink $\phi_1$.
In terms of $u(t,r)$ defined in \eqref{ls}, equation \eqref{eq} takes the form
\be \label{eqw}
\ddot u + L_1 u = f(u,r)\,,
\ee
where the linear operator $L_1$ is defined in \eqref{eq1} and the nonlinear term
\be \label{fwr}
f(u,r)=\sqrt{r^2+a^2} \left[\sin(2\phi_1)-\sin\left(2\phi_1+\frac{2u}{\sqrt{r^2+a^2}}\right)\right] + 2 \cos(2\phi_1) u
\ee
 is real-analytic in $u$ and $r$ (in what follows we suppress the dependence of $f$ on $r$).  The Taylor series of $f(u)$  starts from the quadratic term:
\be\label{tayloru}
f(u)=\frac{2\sin(2\phi_1)}{\sqrt{r^2+a^2}} \, u^2 + \frac{4\cos(2\phi_1)}{2(r^2+a^2)} \, u^3 +\mathcal{O}(u^4).
\ee
As follows from \eqref{expans_inf}, the coefficients in this expansion are decaying functions (exponentially for even powers of $u$ and algebraically for odd powers), which means that the nonlinear terms are spatially localized.  This property, intimately related to the fact that our equation descends from higher dimensions,   will play an important role in our analysis\footnote{We point out that the lack of spatial localization of nonlinear interactions is one of the major difficulties in the studies of asymptotic stability of topological solitons in one spatial dimension, see \cite{KMM, lls, gp}.}.

We recall from the previous section that for $a<a^*$ the spectrum of the operator $L_1$ is purely continuous. In this case, the nonlinear term $f(u)$,  due to its strong spatial localization, does not affect the leading order asymptotic behavior of small amplitude solutions. Consequently, such solutions decay as $t^{-3/2}$ for $t\rightarrow \infty$, which is the three-dimensional  free linear dispersive decay or, equivalently, the one-dimensional linear dispersive decay in the presence of the rapidly decreasing potential (which has no bound states nor a resonance at the bottom of the continuous spectrum) \cite{Schlag, EKMT}.

  For the rest of this section we shall focus on the more interesting case $a>a^*$ where
\[
\mathrm{spec}(L_1)=\{\omega^2\}\cup [2, +\infty),\quad 0<\omega^2<2.
\]
In what follows, the normalized eigenfunction associated to the eigenvalue $\omega^2$ is denoted by $v$, while the modes of the continuous spectrum are denoted by $\eta$.
We decompose the solution as the orthogonal sum of the discrete and continuum modes of~$L_1$
\be \label{decomp}
u(t,r)= \alpha(t) v(r) + \eta(t,r), \quad \mbox{where} \quad \langle v,\eta\rangle=0.
\ee
Substituting this decomposition into  \eqref{eqw} and projecting on the discrete and continuous components, using the projection operators $P f= \langle v,f\rangle v$ and $P^{\perp} f= f-\langle v,f\rangle v $ (where $\langle v,f\rangle :=\int_{-\infty}^{+\infty} \bar v(r) f(r)\,dr$),  we get a system
\begin{eqnarray}
\ddot \alpha +\omega^2 \alpha & = & \langle v, f(\alpha v+\eta)\rangle,\label{eqalpha}\\
\ddot \eta +L_1 \eta &=& P^{\perp} f(\alpha v+\eta) \label{eqeta}\,.
\end{eqnarray}
This system describes  interaction between the internal mode and radiation.
When the nonlinearity is ``switched off", equations \eqref{eqalpha} and \eqref{eqeta} decouple and the internal mode performs harmonic oscillations with frequency $\omega$. If the  perturbation $u$ is small, the system is  weakly coupled and the energy is slowly transferred from the internal mode to the continuum modes and then disperses to infinity. As the result, the amplitude of the internal mode decays asymptotically to zero and the solution converges (on any compact spatial interval) to the static solution $\phi_1$. The key mechanism of this  relaxation process is the nonlinear resonance between the internal mode and the continuum.  It was first  described rigorously by Soffer and Weinstein in a seminal paper \cite{SW}  (see also \cite{PKA} for a formal construction,  \cite{BC,as} for ramifications, and \cite{CM} for a recent review). Below we will adapt their approach to our case, however in contrast to \cite{SW} we will not justify our formal calculations by error estimates (which would be quite technical). Instead,  to feel confident that the results are true, in the next section we will verify them by numerical computations.

 It is convenient to use the complex variable $z=\alpha+\frac{i}{\omega} \dot \alpha$ and rewrite equation \eqref{eqalpha} as the first-order ordinary differential equation
 \be \label{eqzz}
 \dot z + i\omega z = \frac{i}{\omega} \langle v, f\left(\frac{1}{2}(z+\bar z) v+\eta\right)\rangle\,.
 \ee
We expand the nonlinearity in the formal power series in $z, \bar z$ and $\eta$
\be\label{taylor}
f\left(\frac{1}{2}(z+\bar z) v +\eta\right)=
\sum_{k+l+m\geq 2} f_{mkl}\eta^m z^k\bar z^l,
\ee
where  indices run over nonnegative integers. The coefficients $f_{mkl}(r)$ are symmetric in the last two indices and can be read off  from the Taylor series~\eqref{tayloru}.
Substituting this expansion into equations \eqref{eqzz} and \eqref{eqeta} we arrive at the system
\begin{eqnarray}
\dot z +i\omega z &=& \frac{i}{\omega} \sum_{k+l+m\geq 2} \langle f_{mkl}\eta^m, v\rangle\, z^k\bar z^l,
\label{eqzsum}\\
\ddot \eta+L_1\eta &=& \sum_{k+l+m\geq 2} P^\perp (f_{mkl}\eta^m)\,  z^k\bar z^l.
\label{eqetasum}
\end{eqnarray}
In the following we introduce the symbol $\mathcal O_p(z,\eta)$ defined by
\be
\label{op}
\sum_{k+l+m\geq 2} c_{mkl} \eta^m z^k\bar z^l=\mathcal O_p(z,\eta) \quad \mbox{if}\quad \min\{k+l+2m\mid c_{mkl}\neq 0\}=p,
\ee
which incorporates  a heuristic  rule of thumb (to be justified a posteriori) that $\eta=\mathcal{O}(z^2)$.
With this notation, those terms on the right hand sides of equations \eqref{eqzsum} and \eqref{eqetasum}  that do not involve $\eta$ are of the order $\mathcal O_2$, while those that involve $\eta$ are of the order $\mathcal O_3$.

Next, we make a near-identity  coordinate transformation
\be \label{trans1}
\eta = \tilde \eta + \sum_{k+l\geq 2} a_{kl} z^k \bar z^l\,,
\ee
where the coefficients $a_{kl}$ (which are symmetric) are functions of~$r$. The purpose of this transformation is to eliminate the terms of order $\mathcal O_2$ on the right hand side of equation \eqref{eqetasum}. As we will see shortly, this change of variables is formal in the sense that, although a priori $\eta$ is  decaying in space, the functions $a_{kl}$ in general do not decay. To justify rigorously the procedure described below one can consider (\ref{eqw}) in weighted spaces  consistent with \eqref{ls}, c.f. \cite{EKMT} or  in the energy space following the arguments in  \cite{KMM1, KMM4,KMMVDB}. A heuristic reason why the procedure works in the first place is that the system governing the dynamics of the internal modes is always localized by the projection on the eigenspace.

Substituting \eqref{trans1} into \eqref{eqetasum} we get
\be
\ddot{\tilde \eta}+L_1\tilde\eta + \sum_{k+l= 2} (L_1 - (k-l)^2 \omega^2) a_{kl}  z^k\bar z^l=\sum_{k+l= 2} P^\perp f_{0kl} z^k\bar z^l + \mathcal O_3\label{eqetasum2}\,.
\ee
If we manage to remove the $\mathcal O_2$ terms then formally $\tilde \eta=\mathcal O_3$.
To do so we impose the condition
\be \label{eqakl}
\left(L_1-(k-l)^2\omega^2\right) a_{kl}= P^\perp f_{0kl}, \qquad k+l=2.
\ee
For the combination of indices $(k,l)=(1,1)$, there is a unique real-valued solution $a_{11}(r)$ that decays for $|r|\rightarrow \infty$. We do not write it down because the only thing that matters is that this solution is real valued.
For the combination $(k,l)=(2,0)$ we need to consider two cases: (i) $4\omega^2>2$ and (ii) $4\omega^2<2$ (we omit the nongeneric case $4\omega^2=2$ which is more involved).
We postpone the analysis of the case (ii) until afterwards and now we focus on the case (i). In this case, the frequency $2\omega$ lies in the continuous spectrum of $L_1$ and therefore
 solutions of the homogeneous equation $(L_1-4\omega^2) a=0$ are oscillatory at infinity. Among them
there is a unique (complex) solution $a_{20}(r)$ that satisfies the outgoing boundary conditions for $|r|\rightarrow \infty$.  Using the method of variation of parameters,  this solution can be expressed in the form
\be \label{a20}
a_{20}(r)=\frac{i}{2\xi}\, k(r) \int\limits_{-\infty}^r \bar k(s) P^\perp f_{020}(s) ds + \frac{i}{2\xi}\, \bar k(r) \int\limits_r^{\infty} k(s) P^\perp f_{020}(s) ds,
\ee
where $\xi=\sqrt{4\omega^2-2}$ and $k(r)$ is the solution of  the homogeneous equation satisfying the following outgoing boundary condition at $+\infty$
\be \label{jost}
k(r) = e^{i \xi r} m(r),\quad \lim_{r\rightarrow \infty} m(r)=1,
\ee
and $\bar k(r)$ satisfies the corresponding  outgoing boundary condition at $-\infty$.
In \eqref{a20} we used  the Wronskian
\be \label{wron}
W[k(r),\bar k(r)]:=k(r)\bar{k}'(r)-\bar{k}(r) k'(r)=-2i\xi,
\ee
which follows from \eqref{jost}
and the fact that $\bar k(r)=k(-r)$ (what in turn follows from invariance of $L_1$ under reflections $r\rightarrow -r$).

Substituting \eqref{trans1} into \eqref{eqzsum} and using $\tilde \eta=\mathcal O_3$, we get
\be \label{z3}
\dot z +i\omega z = \frac{i}{\omega} \sum_{2\leq k+l\leq 3} \langle v,f_{0kl}\rangle z^k
\bar z^l +  \frac{i}{\omega} \sum_{\substack{k+l=1\\p+q=2}} \langle v,f_{1kl} a_{pq}\rangle z^{k+p}
\bar z^{l+q}+\mathcal O_4\,.
\ee
The key point is that all the dependence on $\tilde \eta$ on  the right hand side is contained in the term $\mathcal O_4$, hence up to the third order this equation is decoupled from the radiation equation.
To factor out fast oscillations with frequency $\omega$ we let $z= e^{-i\omega t} Z$. Substituting this into \eqref{z3}  and dropping all nonresonant terms\footnote{A rigorous justification of this procedure, also called normal form transformation, can be found in \cite{Bam}.} i.e. terms involving  powers of $e^{i\omega t}$ (because such terms time-average to zero), we finally obtain the third order resonant approximation for $t>0$\footnote{We point out that equation \eqref{zres} is not invariant under the time reversal $t\rightarrow -t$.  The arrow of time was selected by  the outgoing boundary conditions imposed on $a_{20}$.}
\be \label{zres}
\dot Z = \frac{i}{\omega} \left(\langle v,f_{021}\rangle+\langle v,f_{110} a_{11}\rangle + \langle v,f_{110} a_{20}\rangle\right) Z^2 \bar Z.
\ee
The first two terms in the bracket  are real, hence multiplying \eqref{zres} by $\bar Z$ and taking the real part we get
\be \label{z2res}
\frac{d}{dt} |Z|^2 = -\Gamma  |Z|^4,\quad \mbox{where}\quad \Gamma= \frac{2}{\omega} \langle v,f_{110} \Im(a_{20})\rangle.
\ee
The coefficient $\Gamma$ can be calculated as follows. From \eqref{tayloru} and \eqref{decomp} we  find
\be \label{ff}
f_{110}= \frac{2 \sin(2\phi_1)}{\sqrt{r^2+a^2}}\, v\quad\mbox{and}\quad f_{020}=\frac{ \sin(2\phi_1)}{2\sqrt{r^2+a^2}}\, v^2=\frac{1}{4} f_{110} v.
\ee
Since  $f_{020}(r)$ is an odd function while $v(r)$ is even, it follows that  $P^{\perp} f_{020}=f_{020}$, and from \eqref{a20} we get
\be \label{a20im}
\Im(a_{20})=\frac{1}{2\xi} \Re\left( k(r) \int\limits_{-\infty}^{\infty} \bar k(s) f_{020}(s) ds\right)\,.
\ee
Inserting this expression  into formula \eqref{z2res} and using \eqref{ff}, we finally obtain
\be\label{Gamma}
\Gamma= \frac{1}{\xi\omega}\, \left\lvert\left\langle k,\frac{\sin(2\phi_1) v^2}{\sqrt{r^2+a^2}}\right\rangle\right\rvert^2\,,
\ee
hence $\Gamma\geq 0$. Generically $\Gamma$ is strictly positive\footnote{This genericity condition is sometimes referred to as the Fermi Golden Rule \cite{S,SW}, which goes back to Dirac's theory of radiation in quantum mechanics \cite{D}.} and  then equation \eqref{z2res} gives
\be \label{decay}
|Z| \sim \Gamma^{-\frac{1}{2}} t^{-\frac{1}{2}}\quad \mbox{as} \quad t\rightarrow \infty.
\ee
The purely imaginary terms on the right hand side of equation \eqref{zres}  determine the phase of $Z$.
Returning to the amplitude $\alpha=\frac{1}{2} \left(e^{-i\omega t} Z +  e^{i\omega t} \bar Z\right)$,
we get the asymptotic behavior
\be \label{final}
\alpha(t) \sim \Gamma^{-\frac{1}{2}} t^{-\frac{1}{2}} \cos\left(\omega t+\theta(t)\right),\quad\mbox{where}\quad \theta(t)=\mathcal{O}(\ln{t}),
\ee
hence as the amplitude of the internal mode decays asymptotically to zero, its frequency tends to the linear frequency $\omega$ (in other words, no frequency shift or `memory' effect occurs).
\vskip 0.2cm
Now, we return to the case (ii) $4\omega^2<2$. In this case the frequency $2\omega$ generated by the quadratic term is below the continuous spectrum so there is no resonant damping present at the third perturbative order (technically, in this case the solution $a_{20}$  of equation \eqref{eqakl} is real and therefore the right hand side of equation \eqref{zres} is purely imaginary).  One needs to go to higher orders to see the damping. Let $N$ be a positive integer such that
\be\label{N}
N^2 \omega^2 < 2 < (N+1)^2 \omega^2\,,
\ee
hence $(N+1)\omega$ is the lowest multiple of the frequency $\omega$ that lies in  the continuous spectrum.
The case (i) discussed above corresponds to $N=1$. For $N\geq 2$ we  iterate the near-identity transformation \eqref{trans1} $N$  times to eliminate terms of order $N+1$ in the radiation equation.
As the result of this iteration,  the internal mode equation decouples from radiation up to  order $\mathcal{O}_{2N+1}$. By the same reasoning as above, this yields the resonant approximation of order $\mathcal{O}_{2N+1}$ for the internal mode equation
\be \label{zresN}
\dot Z = \sum_{1\leq l<N} c_{l} Z^{l+1} \bar{Z}^{l} +  c_N Z^{N+1} \bar{Z}^{N}\,,
\ee
where the coefficients $c_l$ with $l<N$ are purely imaginary, while~$\Re(c_N)~\leq 0$. Thus,
\be\label{ZN}
\frac{d}{dt} |Z|^2 = 2\Re(c_N) |Z|^{2N+2},
\ee
which gives (assuming that $\Re(c_N)$ is strictly negative) for $t\rightarrow \infty$
\be \label{decayN}
|Z| \sim C_N t^{-\frac{1}{2N}},\quad \mbox{where}\quad C_N=\left(2N |\Re(c_N)| \right)^{-\frac{1}{2N}}\,.
\ee
Assuming that the coefficient $c_1$ in \eqref{zresN} is nonzero, the amplitude of the internal mode behaves asymptotically as
\be \label{finalN}
\alpha(t) \sim C_N t^{-\frac{1}{2N}} \cos\left(\omega t+\theta(t)\right),\quad\mbox{where}\quad \theta(t)=\mathcal{O} (t^\frac{N-1}{N}).
\ee
\vskip 0.1cm
Having the formulae for $\alpha(t)$, we now return to  the radiation field $\eta(t,r)$. The asymptotic behavior of $\eta(t,r)$ for $t\rightarrow \infty$ is determined by the quadratic term in equation \eqref{eqetasum}, hence to the leading order we have
\be \label{etaas}
\ddot \eta +L_1 \eta \simeq P^{\perp}\left(\alpha^2\,f_2 v^2+ 2\alpha\, f_2 v \eta+ f_2 \eta^2\right),
\ee
where $f_2=\frac{2 \sin(2\phi_1)}{\sqrt{r^2+a^2}}$ is the coefficient of the quadratic term in \eqref{tayloru}.
All terms on the right hand side are exponentially localized in space and  the dominant contribution comes from the first term which behaves as $\alpha^2$. Inserting $\alpha(t)\sim t^{-1/2N}$ and noting that the solution of the homogeneous equation decays as $t^{-3/2}$ \cite{EKMT}, we conclude that $\eta(t,r)\sim t^{-1/N}$. This heuristic argument justifies a posteriori that $\eta=\mathcal{O}(z^2)$, as claimed above.

\section{Numerical evidence}
In this section we solve equation \eqref{eq}  numerically using  the hyperboloidal formulation of the initial value problem \cite{anil}. As in \cite{BK}
we define  new   coordinates
\begin{equation}\label{hyper_coordinates}
  s=\frac{t}{a}-\sqrt{\frac{r^2}{a^2}+1}\,,\quad y=\arctan\left(\frac{r}{a}\right)\,.
\end{equation}
 The hypersurfaces of constant~$s$ are `hyperboloidal', that is they are spacelike hypersurfaces that approach the `left' and `right' future null infinities of the wormhole spacetime  along the outgoing null cones.  In terms of the coordinates $(s,y)$ and $h(s,y)=\phi(t,r)$ equation \eqref{eq} takes the form
   \be\label{eqy}
    \partial^2_s h + 2 \sin{y} \,\partial_s\partial_y h + \frac{1+\sin^2{y}}{\cos{y}} \,\partial_s h=  \cos^2{\!y} \,\partial^2_y h -a^2\,\frac{\sin(2h)}{\cos^2{y}}\,.\\
    \ee
  We solve this equation for smooth initial data of degree $n$
  \be \label{ic}
  h(0,y)=\alpha(y),\quad \partial_s h(0,y)=\beta(y),
  \ee
   where
 the functions  $\alpha(y)$ and $\beta(y)$ tend exponentially to $\alpha(-\frac{\pi}{2})=0$, $\alpha(\frac{\pi}{2})=n\pi$ and $\beta(\pm\frac{\pi}{2})=0
  $. In this formulation, the $n$-kink denoted by $h_n(y)$ satisfies the boundary conditions $h_n( -\frac{\pi}{2})=0$, $h_n( \frac{\pi}{2})=n\pi$.
 No boundary conditions are imposed because
the principal part of  equation \eqref{eqy} degenerates to $\partial_{s} (\partial_{s}\pm 2\partial_y) h$
at the endpoints $y=\pm \pi/2$, hence there are no ingoing characteristics at the boundaries\footnote{We remark that for massless fields considered in \cite{BK} there was an outflow of energy defined on hyperboloidal slices through the boundaries due the outgoing radiation. In the case at hand, the boundaries do not participate in the evolution, and hence the energy is conserved,  because the group velocity of waves is strictly less than one.}.
\vskip 0.1cm
Following \cite{BK,anil} we define the auxiliary variables
\begin{equation}\label{aux}
  q=\partial_y h \quad\mbox{and}\quad p=\partial_{s} h + \sin{y}\, \partial_y h\,, \nonumber
\end{equation}
  and rewrite equation \eqref{eqy} as the first order symmetric hyperbolic system
  \begin{subequations}\label{symhyp}
  \begin{eqnarray}
  \partial_{s} h &=& p - q\sin{y}\,,\\
  \partial_{s} q &=& \partial_y \left(p - q\sin{y}\right)\,,\\
  \partial_{s} p &=& \partial_y \left(q - p\sin{y}\right) + 2\tan{y} \left(q - p\sin{y}\right)- a^2\,\frac{\sin(2h)}{\cos^2{y}}.
  \end{eqnarray}
\end{subequations}
The initial data \eqref{ic} translate to
\be \label{idata}
  h(0,y)=\alpha(y),\quad q(0,y)=\alpha'(y),\quad p(0,y)=\beta(y)+\sin{y}\, \alpha'(y)\,.
\ee
We solve this system numerically using the method of lines with a fourth-order Runge-Kutta time integration and eighth-order spatial finite differences. One-sided stencils are used at the boundaries. Kreiss-Oliger dissipation  is added in the interior in order to reduce unphysical high-frequency noise.
To suppress   violation of the constraint $q-\partial_y h=0$, we add the term
 $\epsilon (q-\partial_y h)$ with a small negative $\epsilon$ to the right hand side of equation (61b).
\vskip 0.1cm
For any initial data of degree $n$, we find that the solution $h(s,y)$ converges pointwise to the $n$-kink $h_n(y)$ as $s\rightarrow \infty$.
 In the following, we focus on solutions of degree one and  illustrate the results of numerical computations  for sample initial data of the form
 \be\label{sample}
   \alpha(y)= h_1(y)+ e^{-\frac{1}{4}\tan^2{y}},\quad \beta(y)=0.
 \ee
 The convergence rates for several values of $a$ are depicted in Fig.~6.

 \begin{figure}[h]
\includegraphics[width=0.77\textwidth]{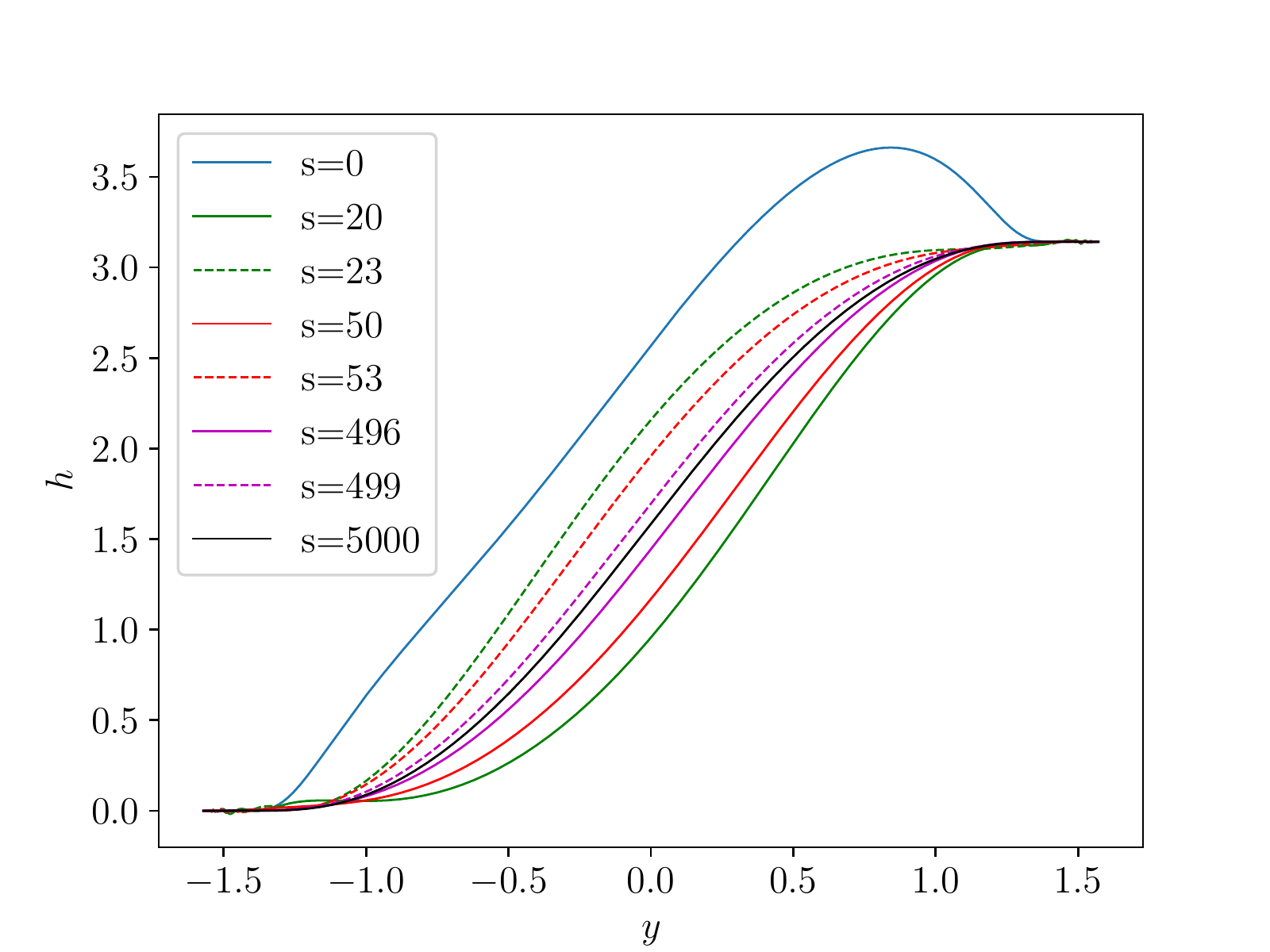}
    \caption{{\small{The snapshots of $h(s,y)$ for the initial data \eqref{sample}. The solution converges to the $1$-kink $h_1(y)$ (solid black line) in an oscillatory manner.}}}
    \end{figure}

 \begin{figure}[h]
\includegraphics[width=0.77\textwidth]{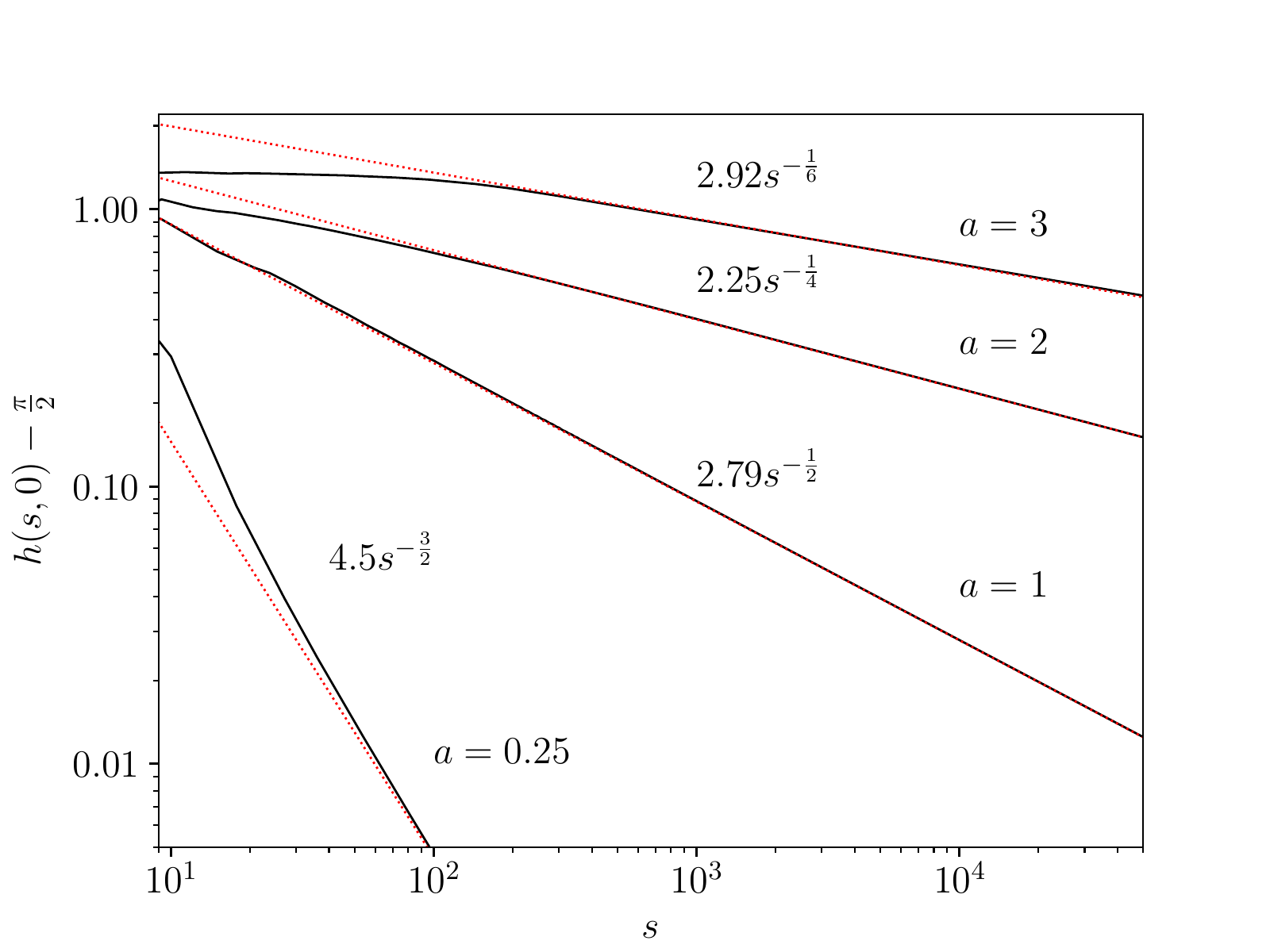}
    \caption{{\small{Amplitudes of perturbations evaluated at $y=0$ in the evolution of initial data \eqref{sample} for different values of the wormhole radius: $a=0.25$ (no gap eigenvalue) and  $a=1$ ($\omega=1.0682$), $a=2$ ($\omega=0.6345$), $a=3$ ($\omega=0.4455$), corresponding to $N=1,2,3$ in \eqref{N}, respectively. The red dashed lines depict analytic predictions with empirically fitted coefficients. }}}
    \end{figure}

The results are in accord with the formulae \eqref{final} and \eqref{finalN} for the decay of the internal mode, and verify  the decay rate $s^{-3/2}$ when the internal mode is absent. We emphasize that both the decay rates and the coefficients are universal (i.e., independent of initial data). In the $N=1$ case the coefficient $2.79$ obtained from the empirical fit agrees (to three decimal places) with  $\Gamma^{-1/2}$ calculated from formula \eqref{Gamma}. This excellent quantitative agreement between analytic and numerical results makes us feel confident that \emph{both} computations are correct.

For the evolutions depicted in Fig.~6 the frequencies of internal modes are well separated from the threshold values $\sqrt{2}/(N+1)$. If the frequency of the internal mode is near a threshold value,   an intermediate dynamics is more complicated. This is illustrated in Fig.~7 where we plot the effective frequency $\omega(s)$ (computed from the distances of subsequent maxima of oscillations) and amplitude of perturbation for $a=1.65$ for which $\omega=0.7422$ is a little above the $N=1$ threshold value $\sqrt{2}/2\approx 0.7071$.
Initially, the effective frequency  is below $\sqrt{2}/2$ and the amplitude decays approximately as $s^{-1/4}$ (as in the $N=2$ case). For later times, the effective frequency increases above the threshold  and concurrently the decay rate of the amplitude undergoes a transition to the asymptotic rate $s^{-1/2}$ ($N=1$ case).

 \begin{figure}[h]
\includegraphics[width=0.49\textwidth]{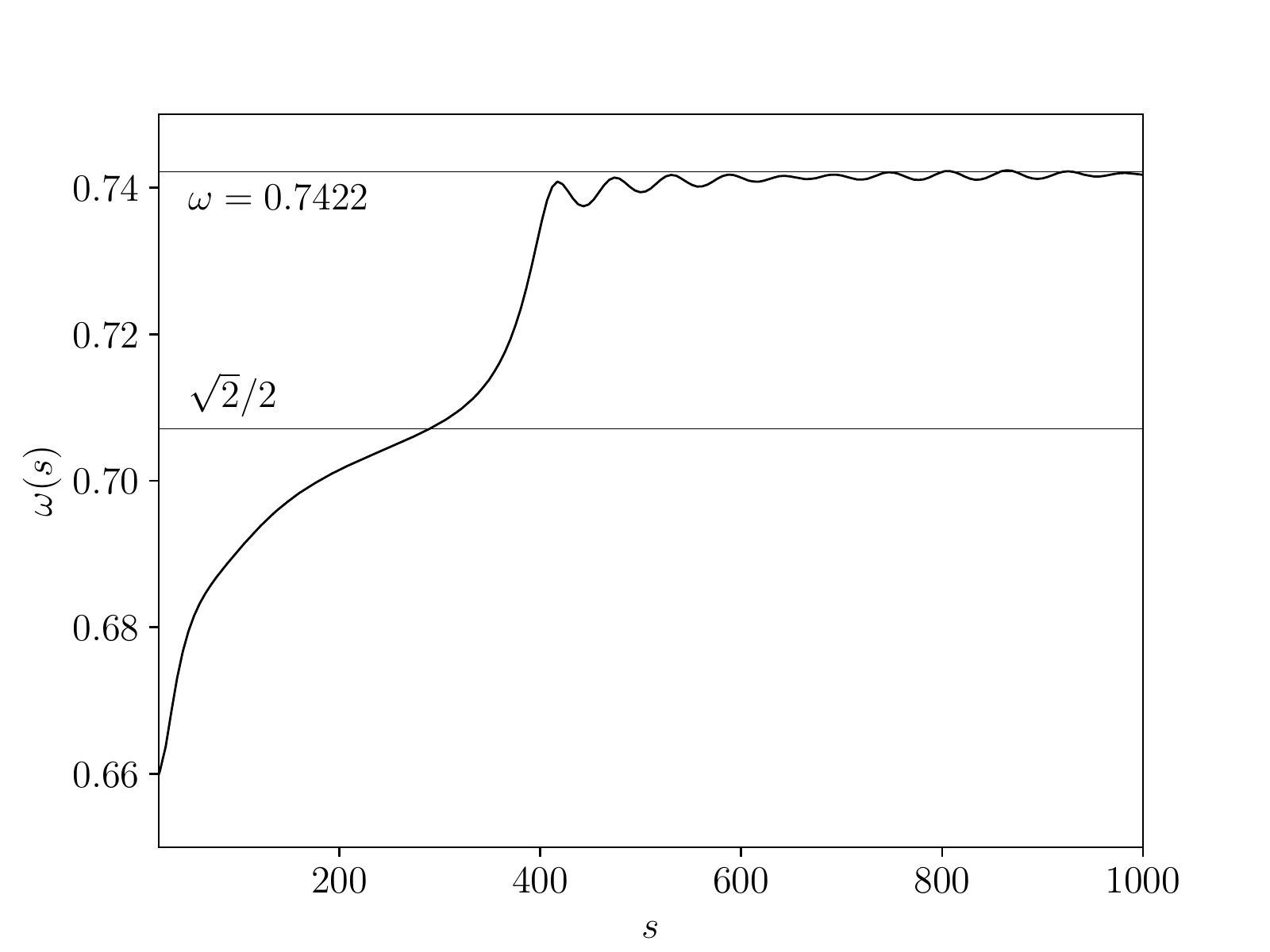}
\includegraphics[width=0.49\textwidth]{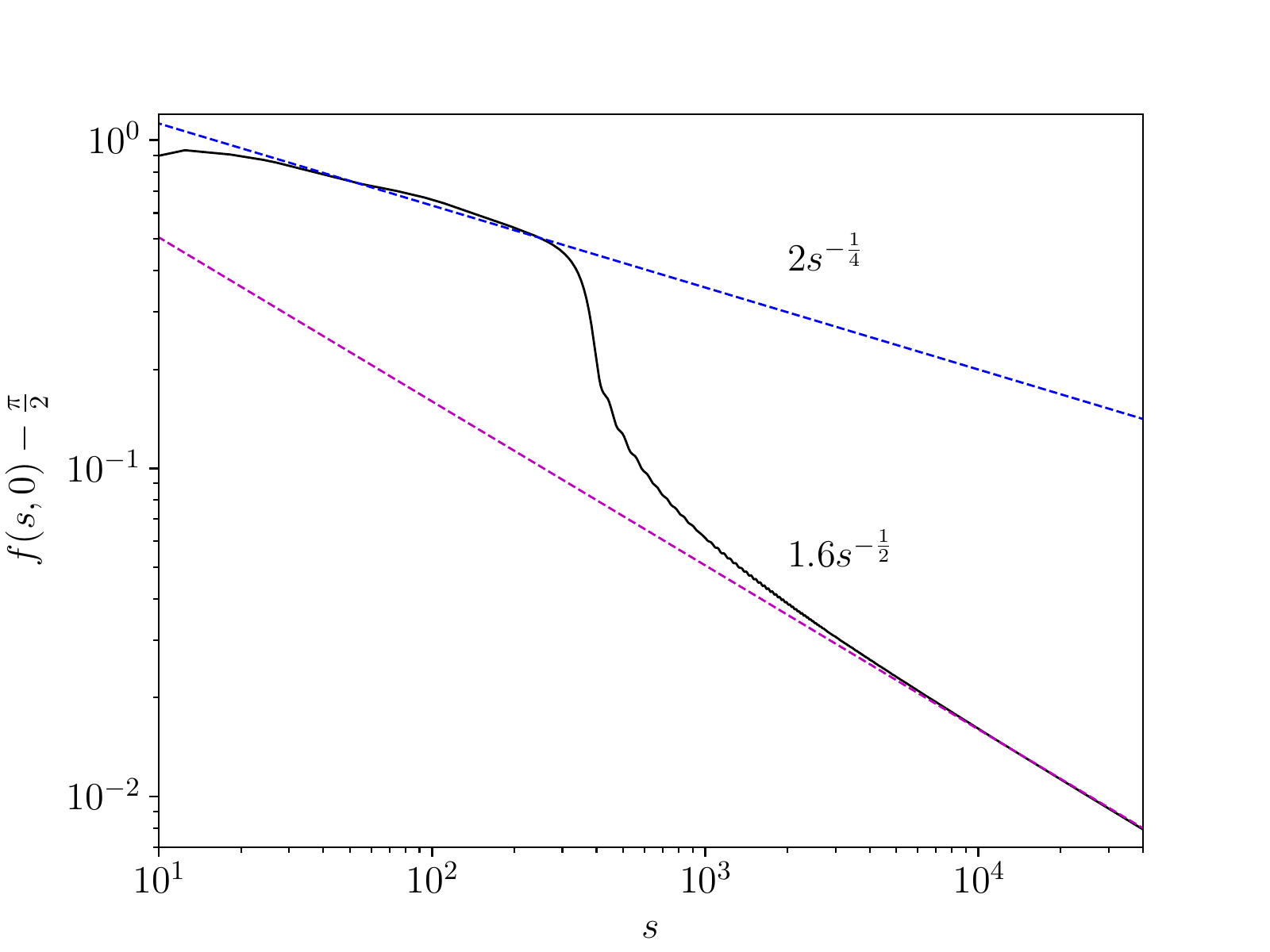}
    \caption{{\small{The effective frequency (left panel) and the amplitude (right panel) in the evolution of initial data \eqref{sample} for $a=1.65$. The late-time small oscillations of the frequency are believed to be due to nonresonant interactions.}}}
    \end{figure}

To summarize, the Sine-Gordon equation on the wormhole is a rich  model for developing understanding of the asymptotic stability of topological solitons with  internal modes. In this work we focused mainly on dynamics of perturbations of the 1-kink. It would be interesting to generalize the weakly nonlinear perturbation analysis from section~4 to $n$-kinks with multiple internal modes. We leave this to future work.

\vskip 0.2cm
\emph{Acknowledgement.} We acknowledge helpful conversations with Gary Gibbons during the early stage of this project.
 This work was supported in part by the National Science Centre grant no.\ 2017/26/A/ST2/00530 (to PB), STFC consolidated grant no. ST/P000681/1 (to MD), FONDECYT 1170164 and CMM Conicyt PIA AFB170001 (to MK).

\section{Appendix}

We describe here how we computed  the coefficients $c_n$ in the asymptotic expansion \eqref{expans_inf}. Since the solutions $\phi_n(r)$ of equation \eqref{ode} are symmetric, it is sufficient to consider the asymptotic behavior at one end, say $r\rightarrow\infty$
\begin{equation}\label{inf2}
  n\pi - \phi_n(r)\sim  \frac{c_n}{r}\,e^{-\sqrt{2} r}.
\end{equation}
Let us observe that for large values of $r$ we have
\begin{equation}\label{lin_as}
  n\pi - \phi_n(r)=c_n \,\phi_L(r)+\mathcal{O}(e^{-3\sqrt{2}r}),
\end{equation}
where $\phi_L(r)$ is the solution of the linearized equation
\begin{equation}\label{ode_lin}
  \phi_L''+\frac{2r}{r^2+a^2}\, \phi_L'- 2\phi_L =0,
\end{equation}
such that \begin{equation}\label{as_lin}
 \phi_L(r)\sim   \frac{1}{r}\,e^{-\sqrt{2} r}\quad  \text{for}\,\, r\rightarrow \infty.
 \end{equation}
Equation \eqref{ode_lin} has two regular singular points $r=\pm a i$ and the irregular singularity at $r=\infty$
and the general solution can be expressed in terms of confluent Heun functions \cite{DLMF}
\begin{equation}\label{heun}
  \phi_L=C_1\,\text{HeunC}\left(0,-\alpha,0,\gamma,\delta,-r^2/a^2\right)    + C_2 \, \text{HeunC}\left(0,\alpha,0,\gamma,\delta,-r^2/a^2\right)\,,
\end{equation}
where $\alpha=\frac{1}{2}, \gamma=\frac{1}{2} a^2,\delta=-\frac{1}{2} a^2+\frac{1}{4}$. In order to find the coefficients $C_1$ and $C_2$ such that $\phi_L$ satisfies \eqref{as_lin} we proceed in two steps. In the first step, we temporarily set $C_1=1$ and determine $C_2$ numerically by solving the equation $\phi_L(r)=0$ for some large value of $r$, say $r=60$. This computation has to be performed with very high precision (up to 60 decimal places) which is possible in Maple where the Heun functions are tabulated. Having that, in the second step, we compute numerically the limit  $\beta =\lim_{r\rightarrow\infty} r e^{\sqrt{2} r} \phi_L(r)$. This is done by evaluating $r_j e^{\sqrt{2} r_j} \phi_L(r_j)$ for an increasing series of $r_j$ (say from $r_1=20$ to $r_{10}=30$) and using the Shanks transformation to accelerate the convergence. Rescaling the coefficients $C_k\mapsto C_k/\beta$ we get the solution $\phi_L(r)$ satisfying \eqref{as_lin}. Finally, taking the numerical solution $\phi_n(r)$ obtained by the shooting method we compute the coefficients $c_n$  using the formula
\begin{equation}\label{cn}
  c_n \approx \frac{n\pi - \phi_n(r_0)}{\phi_L(r_0)},
\end{equation}
where $r_0$ is small enough, say $r_0=6$, so that the value $\phi_n(r_0)$ is accurate.

\end{document}